\documentclass[10pt,twocolumn,a4paper]{article}

\usepackage[left=20mm, right=15mm, top=15mm, bottom=15mm]{geometry}

\usepackage{subfig}
\usepackage{flushend}
\usepackage{indentfirst}
\usepackage{graphics}
\usepackage{amsmath}
\usepackage{graphicx}
\usepackage{epstopdf}
\usepackage{float}

\usepackage[font=footnotesize,labelfont=bf]{caption}
\usepackage[labelsep=period]{caption}

\graphicspath{{./figure/}}

\newtheorem{thm}{Theorem}
\newtheorem{ozn}{Definition}
\begin{document}

\title{\huge \textbf{Modeling of stochastic processes in $L_p(T)$ using orthogonal polynomials}}

\date{}

\twocolumn[
\begin{@twocolumnfalse}
\maketitle
Universal Journal of Applied Mathematics 2(3): 141-147, 2014

DOI: 10.13189/ujam.2014.020304

\vspace{20pt}

\author{\textbf{Oleksandr Mokliachuk}$^{1,*}$\\\\
\footnotesize $^{1}${National Technical University of Ukraine ``Igor Sicorsky Kyiv Polytechnic Institute'', Kyiv, Ukraine}\\
\footnotesize $^{*}$Corresponding email: omoklyachuk@gmail.com}\\\\\\

\end{@twocolumnfalse}
]

\noindent \textbf{\large{Abstract}} \hspace{2pt}
In this paper, models that approximate stochastic processes from the space $Sub_\varphi(\Omega)$ with given reliability and accuracy
in $L_p(T)$ are considered for some specific functions $\varphi(t)$. For processes that are decomposited in series using orthonormal bases,
such models are constructed in the case where elements of such decomposition cannot be found explicitly.
\\

\noindent\textbf{\large{Keywords}} \hspace{2pt}
Models of stochastic processes, $\varphi$-sub-Gaussian processes, orthogonal polynomials
\\


\noindent\hrulefill

\section{\Large{Introduction}}

In different applications of the theory of stochastic processes it is vital to construct the model of the studied process. One way to construct the model of a stochastic process is to represent that process as the infinite series with respect to orthonormal polynomial basis and select the sum of first $N$ elements of this series to be the model.

Consider stochastic process of the second order $X=\{X(t),t\in T\}$, $EX(t)=0$ $\forall t\in T$, and let $B(t,s)=EX(t)\overline{X(s)}$ be the correlation function of this stochastic process $X$. Then the next statement holds true.

\begin{thm} \label{basis-dec-main}\cite{koz-roz-turch} (On decomposition of the stochastic process using an orthonormal basis) Let $X(t)$, $t\in T$ be stochastic process of the second order, $EX(t)=0$ $\forall t\in T$, let $B(t,s)=EX(t)\overline{X(s)}$ be the correlation function of $X$, let $f(t,\lambda)$ be some function from $L_2(\Lambda,\mu)$ space, and let $\{g_k(\lambda), k\in Z\}$ be the orthonormal basis in $L_2(\Lambda,\mu)$ space. Then, correlation function $B(t,s)$ is represented in the form

$$B(t,s) = \int_\Lambda f(t,\lambda) f(s,\lambda) d\mu(\lambda)$$

\noindent if and only if the process can be represented in the form

\begin{equation}
\label{basis}
X(t)=\sum_{k=1}^\infty a_k(t)\xi_k,
\end{equation}
\noindent where

\begin{equation}
\label{basis-f}
a_k(t)=\int_\Lambda f(t,\lambda)\overline{g_k(\lambda)}d\mu(\lambda),
\end{equation}

\noindent $\xi_k$ are centered uncorrelated random variables that satisfy the conditions: $E\xi_k=0$, $E\xi_k\xi_l=\delta_{kl}$, $E\xi_k^2=1$.

\end{thm}

When constructing models of the processes, it is difficult or impossible to find $a_k(t)$ explicitly. In that case, we have to use approximations of these elements. Let us now introduce the model of such process.

\begin{ozn}\label{model-basis-ozn} Let stochastic process $X=\{X(t),t\in T\}$ allow decomposition (\ref{basis}). We will call stochastic process $X_N=\{X_N(t),t\in T\}$ model of the process $X$, if

\begin{equation}
\label{model-basis}
X_N(t) = \sum_{k=1}^N \xi_k\hat{a}_k(t),
\end{equation}

\noindent where $\hat{a}_k(t)$ are approximations of functions $a_k(t)$ in the form (\ref{basis-f}), $\xi_k$ are centered uncorrelated random variables, $E\xi_k=0$, $E\xi_k\xi_l=\delta_{kl}$, $E\xi_k^2=1$.

\end{ozn}

We will consider in details the case where $\xi_k$ are independent $\varphi$-sub-Gaussian random variables.

Since $\hat{a}_k(t)$ are approximations of functions $a_k(t)$, they will introduce some error in the model of stochastic process. The next theorem deals with that case.

\begin{thm}
\label{lp_g2}\cite{mokl}

 Let $X\in S u b_\varphi (\Omega)$, $X=\{X(t),t\in [0,T]\}$ be a stochastic process,

$$\varphi(t)=\left\{\begin{array}{c}\frac{t^2}{\gamma}, t<1\\ \frac{t^\gamma}{\gamma},t\geq1\end{array}\right.,$$

\noindent where  $\gamma>2$. Let  $$c_N=\int_0^T\left(\sum_{k=1}^N \tau^2_\varphi(\xi_k)\delta^2_k(t) +\right.$$ $$ \left.+\sum_{k=N+1}^\infty\tau^2_\varphi(\xi_k)a^2_k(t)  \right)^{p/2} dt<\infty.$$ Model $X_N$ approximates stochastic process
$X$ with given reliability $1-\alpha$ and accuracy $\delta$ in the space
$L_p(0,T)$, if

$$\left\{\begin{array}{c}c_N\leq \delta/(\beta ln \frac{2}{\alpha})^{p/\beta}\\ c_N<\delta/p^{p(1-1/\gamma)}\end{array}\right.,$$

\noindent and $\beta$ satisfies the inequality $1/\beta+1/\gamma=1$.

\end{thm}

A system of orthonormal functions can be considered as basis $\{g_k(\lambda)\}$. It is interesting to consider bases that consist of sets of orthogonal polynomials. The classical examples of such polynomial sets are Chebyshev polynomials, Legendre polynomials, Hermite polynomials, Laguerre polynomials, Jacobi polynomials.

For a system  of polynomials $\{P_n(x)\}$ that are orthogonal with the weight function $h(x)$ on some interval $(a,b)$ and for some fixed $x\in(a,b)$ the next series can be considered:

$$GF(x,\omega)=\sum_{n=0}^\infty \frac{P_n(x)}{n!}\omega^n.$$

Under some minimal conditions, this series has positive convergence radius. In such a case, the function $GF(x,\omega)$ is called generating function of the polynomial set $\{P_n(x)\}$ \cite{suetin}.

If for some functional basis $\{g_k(\lambda)\}$ a generating function exists, namely, if for some $\omega$

$$GF_g(x,\omega)=\sum_{k=0}^\infty g_k(\lambda) \omega^k,$$
then, under additional condition $\tau_\varphi(\xi_k)=\tau \omega^k$, the next equality holds true for the process $X(t)$:

$$\tau_\varphi(X(t)) =\tau_\varphi\left( \sum_{k=0}^\infty \xi_k \int_a^b f(t,\lambda)g_k(\lambda)d\lambda  \right) = $$

$$=\tau_\varphi\left(  \int_a^b f(t,\lambda)\left(\sum_{k=0}^\infty \xi_kg_k(\lambda)\right)d\lambda  \right)\leq$$

$$\leq \tau \int_a^b f(t,\lambda) \left( \sum_{k=0}^\infty\omega^kg_k(\lambda) \right)d\lambda = $$

$$= \tau \int_a^b f(t,\lambda) GF_g(\lambda,\omega) d\lambda.$$

\subsection{Modeling of stochastic processes in $L_p(0,T)$ using the Hermite polynomials}

Let $X=\{X(t),t \in [0,T]\} \in Sub_\varphi(\Omega)$ be stochastic process of the second order, $EX(t)=0$. Let the correlation function of the process $X$, $B(t,s)=EX(t)\overline{X(s)}$, be represented as $$B(t,s)=\int_{-\infty}^\infty f(t,\lambda)f(s,\lambda)d\lambda,$$ where $f(t,\lambda)$, $t\in [0,T]$, $\lambda\in R$ is a family of functions from $L_2(R)$. Since Hermite functions \cite{hermite} form an orthonormal basis, the stochastic process $X$, according to the theorem \ref{basis-dec-main}, can be represented as

$$X(t) = \sum_{k=0}^{\infty}\xi_k \int_{-\infty}^{\infty} f(t,\lambda)\hat{H}_k(\lambda)d\lambda,$$

\noindent where $\xi_k$ are centered uncorrelated random variables, $E\xi_k=0$, $E\xi_k\xi_l=\delta_{kl}$, $E\xi_k^2=1$; $\hat{H}_k(\lambda)$ are Hermite functions:

\begin{equation}
\label{hermite_func}
\hat{H}_k(\lambda)=\frac{H_k(\lambda)}{\sqrt{k!}}\frac{1}{\sqrt[4]{2\pi}}\exp\{-\frac{\lambda^2}{2}\},
\end{equation}

\noindent where $H_k(\lambda)$ are the Hermite polynomials:

$$H_k(\lambda)=(-1)^ke^{\lambda^2/2}\frac{d^n}{d\lambda^n}e^{-\lambda^2/2}.$$


\begin{thm}
\label{hermdec}
 Let a stochastic process  $X=\{X(t),t\in [0,T]\}$ belong to the space $ S u b_\varphi (\Omega)$ with

$$\varphi(t)=\left\{\begin{array}{c}\frac{t^2}{\gamma}, t<1\\ \frac{t^\gamma}{\gamma},t\geq1\end{array}\right.$$

\noindent for  $\gamma>2$, let process $X(t)$ can be represented in the form (\ref{basis}), and let the series $\{\hat{H}_k(t)\}$ of Hermite functions be the basis. Let $$c_N = \int_0^T \left( \int_{-\infty}^\infty Z^2_f(t,\lambda)d\lambda \sum_{k=N+1}^\infty \frac{\tau^2_\varphi(\xi_k)}{k^2+3k+2} +\right.$$ $$\left.+ \sum_{k=1}^N \tau^2_\varphi(\xi_k)\delta^2_k(t)\right)^{p/2}dt<\infty,$$
$$Z_f(t,\lambda) = \left|\frac{\partial^2 f(t,\lambda)}{\partial \lambda^2} - \lambda \frac{\partial f(t,\lambda)}{\partial \lambda}+\frac{\lambda^2-2}{4}f(t,\lambda)\right|,$$
 where function $f(t,s)$ is twice differentiable and bounden with respect to the variable $s$, $Z_f(\lambda)$ is integrable on $R$. The model $X_N(t)$, defined in (\ref{model-basis}), approximates the stochastic process $X(t)$ with given reliability $1-\alpha$ and accuracy $\delta$ in $L_p(0,T)$ spaces, if

$$\left\{\begin{array}{c}c_N\leq \delta/(\beta ln \frac{2}{\alpha})^{p/\beta} \\ c_N<\delta/p^{p(1-1/\gamma)}\end{array}\right.,$$

\noindent where $1/\gamma+1/\beta=1$.

\end{thm}

{\bf Proof. } According to the theorem conditions,

$$a_k(t) = \int_{-\infty}^\infty f(t,\lambda)\hat{H}_k(\lambda)d\lambda = $$$$= \frac{1}{\sqrt[4]{\pi}}\int_{-\infty}^\infty f(t,\lambda)\frac{H_k(\lambda)e^{-\lambda^2/4}}{\sqrt{k!}}d\lambda.$$

\noindent Using properties of Hermite polynomials\cite{hermite},  we can show that

$$\frac{\partial H_k(t)}{\partial t} = k H_{k-1}(t).$$

\noindent Using integration by parts, we get:

$$a_k(t) = \frac{1}{\sqrt[4]{\pi}}\int_{-\infty}^\infty f(t,\lambda)\frac{H_k(\lambda)e^{-\lambda^2/4}}{\sqrt{k!}}d\lambda = $$

$$=\frac{1}{\sqrt[4]{\pi}}\int_{-\infty}^\infty f(t,\lambda)\frac{e^{-\lambda^2/4}}{\sqrt{k+1}\sqrt{(k+1)!}}\frac{\partial H_{k+1}(\lambda)}{\partial \lambda}d\lambda = $$

$$ = \left.\frac{1}{\sqrt[4]{\pi}} f(t,\lambda)\frac{e^{-\lambda^2/4}}{\sqrt{k+1}\sqrt{(k+1)!}}H_{k+1}(\lambda)\right|_{\lambda=-\infty}^{\lambda=\infty} -$$

$$- \frac{1}{\sqrt[4]{\pi}}\int_{-\infty}^\infty \frac{\partial (f(t,\lambda) e^{-\lambda^2/4})}{\partial \lambda}\frac{H_{k+1}(\lambda)}{\sqrt{k+1}\sqrt{(k+1)!}}d\lambda.$$

Since $H_k(\lambda)\exp\{-\lambda^2/4\}$ tends to zero as $\lambda\to\pm\infty$, and $f(t,\lambda)$ is bounded, we get

$$a_k(t)=- \frac{1}{\sqrt[4]{\pi}}\int_{-\infty}^\infty \frac{\partial (f(t,\lambda) e^{-\lambda^2/4})}{\partial \lambda}\frac{H_{k+1}(\lambda)}{\sqrt{k+1}\sqrt{(k+1)!}}d\lambda.$$

\noindent Integration by parts one more time gives

$$a_k(t)=- \frac{1}{\sqrt[4]{\pi}}\int_{-\infty}^\infty \frac{\partial( f(t,\lambda) e^{-\lambda^2/4})}{\partial \lambda}\frac{H_{k+1}(\lambda)}{\sqrt{k+1}\sqrt{(k+1)!}}d\lambda=$$

$$=-\frac{1}{\sqrt[4]{\pi}}\int_{-\infty}^\infty \frac{\partial (f(t,\lambda) e^{-\lambda^2/4})}{\partial \lambda}\times$$$$\times\frac{1}{\sqrt{(k+1)(k+2)}\sqrt{(k+2)!}}\frac{\partial H_{k+2}(\lambda)}{\partial \lambda}d\lambda = $$

$$ = -\left.\frac{1}{\sqrt[4]{\pi}} \frac{\partial (f(t,\lambda) e^{-\lambda^2/4})}{\partial \lambda}\times\right.$$$$\left.\times\frac{1}{\sqrt{(k+1)(k+2)}\sqrt{(k+2)!}}H_{k+1}(\lambda)\right|_{\lambda=-\infty}^{\lambda=\infty} + $$

$$+ \frac{1}{\sqrt[4]{\pi}}\int_{-\infty}^\infty \frac{\partial^2 (f(t,\lambda) e^{-\lambda^2/4})}{\partial \lambda^2}\times$$$$\times\frac{H_{k+2}(\lambda)}{\sqrt{(k+1)(k+2)}\sqrt{(k+2)!}}d\lambda. $$

It is certain that $H_{k+1}(\lambda)\partial (f(t,\lambda)e^{-\lambda^2/4})/\partial \lambda$ tends to zero as $\lambda\to\pm\infty$, since $\partial (f(t,\lambda)e^{-\lambda^2/4})/\partial \lambda $ $=e^{-\lambda^2/4}(\partial f(t,\lambda)/\partial\lambda-\lambda f(t,\lambda))$, $H_{k+1}(\lambda)e^{-\lambda^2/4}\to 0$, and   $\partial f(t,\lambda)/\partial\lambda-\lambda f(t,\lambda)$ is bounded, because $f(t,\lambda)$ and $\lambda f(t,\lambda)$ are bounded due to the conditions of the theorem. Then,

$$a_k(t) = \frac{1}{\sqrt[4]{\pi}}\int_{-\infty}^\infty \frac{\partial^2 (f(t,\lambda) e^{-\lambda^2/4})}{\partial \lambda^2}\times$$$$\times\frac{H_{k+2}(\lambda)}{\sqrt{(k+1)(k+2)}\sqrt{(k+2)!}}d\lambda.$$

\noindent Because $\hat{H}_k$ is an orthonormal basis, $\int_{-\infty}^\infty \hat{H}^2_k(t)dt = 1$. That's why

$$\frac{1}{\sqrt[4]{\pi}}\int_{-\infty}^\infty \frac{\partial^2 (f(t,\lambda) e^{-\lambda^2/4})}{\partial \lambda^2}\times$$$$\times\frac{H_{k+2}(\lambda)}{\sqrt{(k+1)(k+2)}\sqrt{(k+2)!}}d\lambda = $$

$$=\frac{1}{\sqrt{(k+1)(k+2)}}\int_{-\infty}^\infty \frac{\partial^2 (f(t,\lambda) e^{-\lambda^2/4})}{\partial \lambda^2}\times$$$$\times\frac{H_{k+2}(\lambda)e^{-\lambda^2/4}}{\sqrt{(k+2)!}\sqrt[4]{2\pi}}e^{\lambda^2/4}d\lambda=$$

$$ = \frac{1}{\sqrt{(k+1)(k+2)}}\int_{-\infty}^\infty \frac{\partial^2 (f(t,\lambda) e^{-\lambda^2/4})}{\partial \lambda^2}\times$$$$\times g_{k+2}(\lambda)e^{\lambda^2/4}d\lambda\leq$$

$$\leq \frac{1}{\sqrt{(k+1)(k+2)}}\left(\int_{-\infty}^\infty \left(\frac{\partial^2 (f(t,\lambda) e^{-\lambda^2/4})}{\partial \lambda^2}\right.\right.\times$$$$\left.\left.\times e^{\lambda^2/4}\right)^2d\lambda\right)^\frac{1}{2}\left(\int_{-\infty}^\infty g^2_{k+2}(\lambda) d\lambda\right)^\frac{1}{2}=$$

$$= \frac{1}{\sqrt{(k+1)(k+2)}}\left(\int_{-\infty}^\infty \left(\frac{\partial^2 (f(t,\lambda) e^{-\lambda^2/4})}{\partial \lambda^2}\right.\right.\times$$$$\left.\left.\times e^{\lambda^2/4}\right)^2d\lambda\right)^\frac{1}{2}$$

Besides,

$$\frac{\partial^2 (f(t,\lambda) e^{-\lambda^2/4})}{\partial \lambda^2}e^{\lambda^2/4} = \frac{\partial}{\partial \lambda}\left(\frac{\partial f(t,\lambda)}{\partial \lambda}e^{-\lambda^2/4}\right.-$$

$$\left.-\frac{1}{2}e^{-\lambda^2/4}\lambda f(t,\lambda)\right)e^{\lambda^2/4}=\left(\frac{\partial}{\partial \lambda}\left(\frac{\partial f(t,\lambda)}{\partial \lambda}e^{-\lambda^2/4}\right)-\right.$$

$$\left.-\frac{\partial}{\partial \lambda}\left(\frac{1}{2}e^{-\lambda^2/4}\lambda f(t,\lambda)\right)\right)e^{\lambda^2/4}=\left(\frac{\partial^2 f(t,\lambda)}{\partial \lambda^2}e^{-\lambda^2/4}-\right.$$

$$\left.-\frac{1}{2}e^{-\lambda^2/4}\lambda \frac{\partial f(t,\lambda)}{\partial \lambda} - \frac{1}{2}e^{-\lambda^2/4}\lambda \frac{\partial f(t,\lambda)}{\partial \lambda} -\right.$$

$$\left.-f(t,\lambda)\frac{\partial}{\partial \lambda}\left(\frac{1}{2}e^{-\lambda^2/4}\lambda\right)\right)e^{\lambda^2/4} = \left(\frac{\partial^2 f(t,\lambda)}{\partial \lambda^2}e^{-\lambda^2/4}-\right.$$$$\left.-e^{-\lambda^2/4}\lambda \frac{\partial f(t,\lambda)}{\partial \lambda}\right.-$$

$$\left. - f(t,\lambda)\left(\frac{1}{2}e^{-\lambda^2/4}-\frac{1}{4}e^{-\lambda^2/4}\lambda^2\right) \right)e^{\lambda^2/4} = $$

$$=\frac{\partial^2 f(t,\lambda)}{\partial \lambda^2} - \lambda \frac{\partial f(t,\lambda)}{\partial \lambda}+\frac{\lambda^2-2}{4}f(t,\lambda).$$

 Finally, we obtain

$$a_k(t)\leq \left(\frac{1}{(k+1)(k+2)}\int_{-\infty}^\infty Z_f^2(t,\lambda) d\lambda\right)^\frac{1}{2},$$

\noindent where

$$Z_f(t,\lambda) = \left|\frac{\partial^2 f(t,\lambda)}{\partial \lambda^2} - \lambda \frac{\partial f(t,\lambda)}{\partial \lambda}+\frac{\lambda^2-2}{4}f(t,\lambda)\right|.$$

\noindent Therefore,

$$c_N = \int_0^T\left(\sum_{k=1}^N \tau^2_\varphi(\xi_k)\delta^2_k(t) +\right.$$$$\left.+ \sum_{k=N+1}^\infty\tau^2_\varphi(\xi_k)a^2_k(t)  \right)^{p/2} dt \leq $$

$$\leq\int_0^T \left( \int_{-\infty}^\infty Z^2_f(t,\lambda)d\lambda \sum_{k=N+1}^\infty \frac{\tau^2_\varphi(\xi_k)}{k^2+3k+2} \right.+$$$$+\left. \sum_{k=1}^N \tau^2_\varphi(\xi_k)\delta^2_k(t)\right)^{p/2}dt.$$

\noindent Finally, the statement of this theorem is derived from the theorem \ref{lp_g2}.

\begin{thm}

Let stochastic process $X=\{X(t),t\in [0,T]\}$ belong to the space $ S u b_\varphi (\Omega)$,

$$\varphi(t)=\left\{\begin{array}{c}\frac{t^2}{\gamma}, t<1\\ \frac{t^\gamma}{\gamma},t\geq1\end{array}\right.$$

\noindent for  $\gamma>2$, let the process $X(t)$ allow representation in the form (\ref{basis}), and the series $\hat{H}_k(t)$ of Hermite functions is the basis. Let $\tau_\varphi(\xi_k)=\tau \omega^k$, $|\omega|<1$, and

$$c_N = \int_0^T \left( \frac{\tau}{\sqrt{(1-\omega^2)}}\left(\int_{-\infty}^\infty f^2(t,\lambda)d\lambda\right)^{1/2}-\right.$$$$\left.-\sum_{k=0}^N \tau \omega^k\hat{a}_k(t)\right)^{p}dt<\infty.$$

Model $X_N(t)$, provided in (\ref{model-basis}), approximates $X(t)$ with given reliability $1-\alpha$ and accuracy $\delta$ in the space $L_p(0,T)$, if

$$\left\{\begin{array}{c}c_N\leq \delta/(\beta ln \frac{2}{\alpha})^{p/\beta} \\ c_N<\delta/p^{p(1-1/\gamma)}\end{array}\right.,$$

\noindent where $1/\gamma+1/\beta=1$.

\end{thm}

{\bf Proof. }According to the statement of the theorem,

$$a_k(t) = \int_{-\infty}^\infty f(t,\lambda)\hat{H}_k(\lambda)d\lambda = $$ $$=\frac{1}{\sqrt[4]{\pi}}\int_{-\infty}^\infty f(t,\lambda)\frac{H_k(\lambda)e^{-\lambda^2/2}}{\sqrt{k!}}d\lambda.$$

Using the properties of Hermite polynomials $\{H_k(\lambda)\}$, we get

$$GF_{H^2}(\lambda,\omega) = \sum_{k=1}^\infty \frac{H^2_k(\lambda)}{k!2^k}\omega^k = \frac{1}{\sqrt{1-\omega^2}}\exp\left\{\frac{2\lambda^2\omega}{1+\omega}\right\}.$$

\noindent Moreover, such series converges for $|\omega|<1$. Under the conditions of the theorem, $\tau_\varphi(\xi_k)=\tau \omega^k$. That's why, for the process $X(t)$ the next condition holds true:

$$\tau_\varphi(X(t)) = \tau_\varphi\left(\sum_{k=1}^\infty \xi_k\int_{-\infty}^\infty f(t,\lambda)\hat{H}_k(\lambda)d\lambda\right) \leq  $$$$\leq\tau\int_{-\infty}^\infty f(t,\lambda) \sum_{k=0}^\infty \omega^k\hat{H}_k(\lambda)d\lambda\leq$$

 $$\leq \tau \left(\int_{-\infty}^\infty f^2(t,\lambda)d\lambda\right)^{1/2} \times$$$$\times \left(\int_{-\infty}^\infty \left(\sum_{k=0}^\infty\frac{H_k(\lambda)\exp\{-\lambda^2/2\}}{\sqrt{k!2^k\sqrt{\pi}}}\omega^k\right)^2d\lambda\right)^{1/2} \leq $$

 $$\leq \tau \left(\int_{-\infty}^\infty f^2(t,\lambda)d\lambda\right)^{1/2} \times$$$$\times \left(\int_{-\infty}^\infty \sum_{k=0}^\infty\frac{H^2_k(\lambda)\exp\{-\lambda^2\}}{k!2^k\sqrt{\pi}}\omega^{2k}d\lambda\right)^{1/2} =$$

 $$= \tau \left(\int_{-\infty}^\infty f^2(t,\lambda)d\lambda\right)^{1/2} \times$$$$\times \left(\frac{1}{\sqrt{\pi}}\int_{-\infty}^\infty GF_{H^2}(\lambda,\omega^2)\exp\{-\lambda^2\}d\lambda\right)^{1/2} =$$

 $$= \tau \left(\int_{-\infty}^\infty f^2(t,\lambda)d\lambda\right)^{1/2} \times$$$$\times \left(\frac{1}{\sqrt{\pi(1-\omega^4)}}\int_{-\infty}^\infty\exp\left\{\frac{2\lambda^2\omega}{1+\omega}\right\}\exp\{-\lambda^2\}d\lambda\right)^{1/2} = $$

 $$ = \tau \left(\int_{-\infty}^\infty f^2(t,\lambda)d\lambda\right)^{1/2} \times$$$$\times \left(\frac{1}{\sqrt{\pi(1-\omega^4)}}\int_{-\infty}^\infty\exp\left\{\lambda^2\frac{\omega^2-1}{\omega^2+1}\right\}d\lambda\right)^{1/2} = $$

 $$ = \tau \left(\int_{-\infty}^\infty f^2(t,\lambda)d\lambda\right)^{1/2} \times$$$$\times \left(\frac{1}{\sqrt{\pi(1-\omega^4)}} \sqrt{\pi\frac{1+\omega^2}{1-\omega^2}}\right)^{1/2} =$$

 $$ = \frac{\tau}{\sqrt{(1-\omega^2)}}\left(\int_{-\infty}^\infty f^2(t,\lambda)d\lambda\right)^{1/2}.$$

 Given these considerations, the estimator of the model of the process will take the next form:

 $$\tau_\varphi(X(t)-X_N(t)) = \tau_\varphi\left(\sum_{k=0}^\infty \xi_ka_k(t)-\sum_{k=0}^N \xi_k\hat{a}_k(t)\right)= $$$$=\tau_\varphi\left(\sum_{k=0}^N \xi_k\delta_k(t) + \sum_{k=N+1}^\infty \xi_ka_k(t)\right)\leq$$

 $$\leq \sum_{k=0}^N \tau_\varphi(\xi_k)\delta_k(t) + \sum_{k=N+1}^\infty \tau_\varphi(\xi_k) a_k(t)  =$$$$= \sum_{k=0}^\infty \tau \omega^k a_k(t)-\sum_{k=0}^N \tau \omega^k\hat{a}_k(t)\leq$$

 $$\leq \frac{\tau}{\sqrt{(1-\omega^2)}}\left(\int_{-\infty}^\infty f^2(t,\lambda)d\lambda\right)^{1/2}-\sum_{k=0}^N \tau \omega^k\hat{a}_k(t).$$

 Statement of this theorem follows from the theorem \ref{lp_g2} and the last inequality.

\subsection{Modeling of stochastic processes in $L_p(0,T)$ using the Chebyshev polynomials}

Let the process $X(t)$ has the same properties as the process from the previous section. Let orthonotmal Chebyshev polynomials be used as the basis:

$$\hat{T}_n(\lambda) = \sqrt{\frac{2}{\pi}}T_n(\lambda),$$ where

$$T_n(\lambda) = \cos(n \arccos \lambda).$$

In such a case we can proof the next theorem.

\begin{thm}

 Let a stochastic process $X=\{X(t),t\in [0,T]\}$ belong to the space $ S u b_\varphi (\Omega)$ with

$$\varphi(t)=\left\{\begin{array}{c}\frac{t^2}{\gamma}, t<1\\ \frac{t^\gamma}{\gamma},t\geq1\end{array}\right.$$

\noindent for  $\gamma>2$, let $X(t)$ is representated in the form (\ref{basis}), and let the system of orthonormal Chebyshev polynomials $\{\hat{T}_k(t)\}$ be used as the basis. Let $\tau_\varphi(\xi_k)=\tau \omega^k$, $0<\omega<1$, and

$$c_N = \int_0^T \left(  \sqrt{\frac{2}{\pi}}\tau\left(\int_{-1}^1 f^2(t,\lambda)d\lambda\right)^{1/2} \sqrt{D_T(\omega)}-\right.$$$$\left.-\sum_{k=0}^N \tau \omega^k\hat{a}_k(t)\right)^{p}dt<\infty,$$

$$D_T(\omega)=2\frac{1}{{\omega(4+3\omega^2+\omega^4)}}\left(\omega(5+5\omega^2+2\omega^2)+\right.$$

$$+\left.(4+7\omega^2+4\omega^4+\omega^6)\ln\{(\omega^2-\omega+2)/(\omega^2+\omega+2)\}\right).$$

Model $X_N(t)$, determined in (\ref{model-basis}), approximates the process $X(t)$ with given reliability $1-\alpha$ and accuracy $\delta$ in the space $L_p(0,T)$, if

$$\left\{\begin{array}{c}c_N\leq \delta/(\beta ln \frac{2}{\alpha})^{p/\beta} \\ c_N<\delta/p^{p(1-1/\gamma)}\end{array}\right.,$$

\noindent where $1/\gamma+1/\beta=1$.
\end{thm}

{\bf Proof. } Under the conditions of the theorem,

$$a_k(t) = \int_{-1}^1 f(t,\lambda)\hat{T}_n(\lambda)d\lambda = \int_{-1}^1 f(t,\lambda)\sqrt{\frac{2}{\pi}}T_n(\lambda)d\lambda.$$

The generating function of orthogonal Chebyshev polynomials of the first kind $\{T_k(\lambda)\}$ has the next form:

$$GF_T(\lambda,\omega) = \sum_{k=0}^\infty T_k(\lambda)\omega^k = \frac{1-\omega \lambda}{2-\omega\lambda+\omega^2}$$ for $0<\omega<1$. Under the conditions of the theorem, $\tau_\varphi(\xi_k)=\tau \omega^k$. Thats why for the process $X(t)$ the next  condition is true:

$$\tau_\varphi(X(t)) = \tau_\varphi\left(\sum_{k=1}^\infty \xi_k\int_{-1}^1 f(t,\lambda)\hat{T}_k(\lambda)d\lambda\right) \leq$$$$\leq  \tau\int_{-1}^1 f(t,\lambda) \sum_{k=0}^\infty \omega^k\hat{T}_k(\lambda)d\lambda\leq$$

 $$\leq \tau \left(\int_{-1}^1 f^2(t,\lambda)d\lambda\right)^{1/2}  \times$$$$\times \left(\int_{-1}^1 \left(\sum_{k=0}^\infty\sqrt{\frac{2}{\pi}}T_k(\lambda)\right)^2d\lambda\right)^{1/2} = $$

$$ =  \sqrt{\frac{2}{\pi}}\tau\left(\int_{-1}^1 f^2(t,\lambda)d\lambda\right)^{1/2}  \times$$$$\times \left(\int_{-1}^1\left(\frac{1-\omega \lambda}{2-\omega\lambda+\omega^2} \right)^2d\lambda\right)^{1/2}.$$

Lets calculate the second integral of the last expression separetely.

$$\int_{-1}^1\left(\frac{1-\omega \lambda}{2-\omega\lambda+\omega^2} \right)^2d\lambda = $$ $$=2\frac{1}{{\omega(4+3\omega^2+\omega^4)}}\left(\omega(5+5\omega^2+2\omega^2)+\right.$$

$$+\left.(4+7\omega^2+4\omega^4+\omega^6)\ln\{(\omega^2-\omega+2)/(\omega^2+\right.$$$$\left.+\omega+2)\}\right):=D_T(\omega).$$

\noindent Then the estimator of $\tau_\varphi(X(t))$ will take the form:

$$\tau_\varphi(X(t))\leq\sqrt{\frac{2}{\pi}}\tau\left(\int_{-1}^1 f^2(t,\lambda)d\lambda\right)^{1/2} \sqrt{D_T(\omega)}.$$

Given these considerations, the estimator of the model of the process will take the following form:

 $$\tau_\varphi(X(t)-X_N(t)) = \tau_\varphi\left(\sum_{k=0}^\infty \xi_ka_k(t)-\sum_{k=0}^N \xi_k\hat{a}_k(t)\right)=$$$$= \tau_\varphi\left(\sum_{k=0}^N \xi_k\delta_k(t) + \sum_{k=N+1}^\infty \xi_ka_k(t)\right)\leq$$

 $$\leq \sum_{k=0}^N \tau_\varphi(\xi_k)\delta_k(t) + \sum_{k=N+1}^\infty \tau_\varphi(\xi_k) a_k(t)  = $$$$\sum_{k=0}^\infty \tau \omega^k a_k(t)-\sum_{k=0}^N \tau \omega^k\hat{a}_k(t)\leq$$

 $$\leq\sqrt{\frac{2}{\pi}}\tau\left(\int_{-1}^1 f^2(t,\lambda)d\lambda\right)^{1/2} \sqrt{D_T(\omega)}-\sum_{k=0}^N \tau \omega^k\hat{a}_k(t).$$

 Statement of this theorem follows from the theorem \ref{lp_g2} and the last inequality.

We can also proof the similar theorem in the case of Chebyshev polynomials of the second kind:

$$U_n(\lambda)=\frac{\sin((n+1)\arccos \lambda)}{\sqrt{1-\lambda^2}}$$

$$\hat{U}_n(\lambda) = \sqrt{\frac{2}{\pi}}U_n(\lambda).$$ Generating function of the series $\{U_n(\lambda)\}$ has the form

$$GF_U(\lambda,\omega) = \sum_{k=0}^\infty \omega^k U_n(\lambda).$$

Let the process $X(t)$ has the same properties as in the previous theorem Then, the next proposition holds true.

\begin{thm}  Let stochastic process $X=\{X(t),t\in [0,T]\}$ belong to the space $ S u b_\varphi (\Omega)$,

$$\varphi(t)=\left\{\begin{array}{c}\frac{t^2}{\gamma}, t<1\\ \frac{t^\gamma}{\gamma},t\geq1\end{array}\right.$$

\noindent for $\gamma>2$, let process $X(t)$ can be represented in the form (\ref{basis}), and let series $\{\hat{U}_k(t)\}$ of orthonormal Chebyshev polynomials be the basis. Let also $\tau_\varphi(\xi_k)=\tau \omega^k$, $0<\omega<1$, and

$$c_N = \int_0^T \left(  \frac{2\tau}{\sqrt{\pi}(\omega^2-1)}\left(\int_{-1}^1 f^2(t,\lambda)d\lambda\right)^{1/2}-\right.$$$$\left.-\sum_{k=0}^N \tau \omega^k\hat{a}_k(t)\right)^{p}dt<\infty.$$

Model $X_N(t)$, determined in (\ref{model-basis}), approximates the process $X(t)$ with given reliability $1-\alpha$ and accuracy $\delta$ in the space $L_p(0,T)$, if

$$\left\{\begin{array}{c}c_N\leq \delta/(\beta ln \frac{2}{\alpha})^{p/\beta} \\ c_N<\delta/p^{p(1-1/\gamma)}\end{array}\right.,$$

\noindent where $1/\gamma+1/\beta=1$.
\end{thm}

{\bf Proof. } Under the conditions of the theorem, we have:

$$a_k(t) = \int_{-1}^1 f(t,\lambda)\hat{U}_n(\lambda)d\lambda = \int_{-1}^1 f(t,\lambda)\sqrt{\frac{2}{\pi}}U_n(\lambda)d\lambda.$$

Generation function of the series of orthogonal Chebyshev polynomials of the second kind $\{U_k(\lambda)\}$ has the following form:

$$GF_T(\lambda,\omega) = \sum_{k=0}^\infty U_k(\lambda)\omega^k = \frac{1}{1-2\omega\lambda+\omega^2}$$ for $0<\omega<1$. As stated in the theorem's conditions, $\tau_\varphi(\xi_k)=\tau \omega^k$. That's why for stochastic process $X(t)$ the next statement holds true:

$$\tau_\varphi(X(t)) = \tau_\varphi\left(\sum_{k=1}^\infty \xi_k\int_{-1}^1 f(t,\lambda)\hat{U}_k(\lambda)d\lambda\right) \leq$$$$\leq  \tau\int_{-1}^1 f(t,\lambda) \sum_{k=0}^\infty \omega^k\hat{U}_k(\lambda)d\lambda\leq$$

 $$\leq \tau \left(\int_{-1}^1 f^2(t,\lambda)d\lambda\right)^{1/2}\times$$$$\times \left(\int_{-1}^1 \left(\sum_{k=0}^\infty\sqrt{\frac{2}{\pi}}U_k(\lambda)\right)^2d\lambda\right)^{1/2} = $$

$$ =  \sqrt{\frac{2}{\pi}}\tau\left(\int_{-1}^1 f^2(t,\lambda)d\lambda\right)^{1/2} \times$$$$\times \left(\int_{-1}^1\left(\frac{1}{1-2\omega\lambda+\omega^2} \right)^2d\lambda\right)^{1/2}=$$

$$=\sqrt{\frac{2}{\pi}}\tau\left(\int_{-1}^1 f^2(t,\lambda)d\lambda\right)^{1/2} \frac{\sqrt{2}}{(\omega^2-1)} = $$$$=\frac{2\tau}{\sqrt{\pi}(\omega^2-1)}\left(\int_{-1}^1 f^2(t,\lambda)d\lambda\right)^{1/2}.$$

Taking into account these considerations, the estimator of the model of the stochastic process will fit the next condition:

 $$\tau_\varphi(X(t)-X_N(t)) = \tau_\varphi\left(\sum_{k=0}^\infty \xi_ka_k(t)-\sum_{k=0}^N \xi_k\hat{a}_k(t)\right)= $$$$=\tau_\varphi\left(\sum_{k=0}^N \xi_k\delta_k(t) + \sum_{k=N+1}^\infty \xi_ka_k(t)\right)\leq$$

 $$\leq \sum_{k=0}^N \tau_\varphi(\xi_k)\delta_k(t) + \sum_{k=N+1}^\infty \tau_\varphi(\xi_k) a_k(t)  = $$$$=\sum_{k=0}^\infty \tau \omega^k a_k(t)-\sum_{k=0}^N \tau \omega^k\hat{a}_k(t)\leq$$

 $$\leq \frac{2\tau}{\sqrt{\pi}(\omega^2-1)}\left(\int_{-1}^1 f^2(t,\lambda)d\lambda\right)^{1/2}-\sum_{k=0}^N \tau \omega^k\hat{a}_k(t).$$

Statement of this theorem follows from the theorem \ref{lp_g2} and the last inequality.

\section{\Large{Conclusions}}

Theorems are proved that allow to construct models of stochastic processes from $Sub_\varphi(\Omega)$ in the case where these processes are represented in series with respect to some orthogonal polynomial systems. The case when functional components of such decompositions cannot be found explicitly is studied.

\noindent\hrulefill

\end{document}